\documentclass[reqno,b5paper,12pt]{amsart}
\usepackage{amsmath,amssymb,amsthm,mathptmx}

\newtheorem{theorem}{Theorem}[section]
\newtheorem{corollary}[theorem]{Corollary}

\theoremstyle{definition}

\newtheorem{remark}[theorem]{Remark}

\DeclareMathOperator{\RE}{Re}

\begin{document}

\title[Radii of Starlikeness and Convexity]{Radii of Starlikeness and Convexity of Analytic
Functions\\  Satisfying Certain Coefficient Inequalities}

\newcommand{\acr}{\newline\indent}
\author {V. Ravichandran}
\address{Department of Mathematics   \acr
University of Delhi  \acr Delhi 110 007\acr INDIA}
\email{vravi@maths.du.ac.in}

\address{School of Mathematical Sciences  \acr Universiti Sains
Malaysia  \acr 11800 USM Penang\acr MALAYSIA}
\email{vravi68@gmail.com}

\begin{abstract}For $0\leq \alpha <1$, the sharp radii of
starlikeness and convexity of order $\alpha$ for functions of the
form $f(z)=z+a_2z^2+a_3z^3+\cdots$ whose Taylor coefficients $a_n$
satisfy the conditions  $|a_2|=2b$, $0\leq b\leq 1$, and $|a_n|\leq
n $, $M$ or $M/n$ ($M>0$) for $n\geq 3$ are  obtained. Also a class
of functions related to Carath\'eodory functions is considered.
\end{abstract}



\thanks{The problem was suggested by Prof. S. Ponnusamy in  a lecture
at Universiti Sains Malaysia during his recent visit. The author is
thankful to him for helpful discussion during the preparation of the
manuscript.}

 \maketitle

\section{Introduction}
Let $\mathcal{A}$ be the class of  analytic functions $f$ in the
unit disk $\mathbb{D}=\{z\in\mathbb{C} \colon |z|<1\}$ with Taylor
series expansion $f(z)=z+\sum_2^\infty a_nz^n$. For functions
belonging to the subclass $\mathcal{S}$ of  $\mathcal{A}$ consisting
of univalent functions, it is well-known that $|a_n|\leq n$ for
$n\geq 2$. A function $f$ whose coefficients satisfy the inequality
$|a_n|\leq n$ for $n\geq 2$ are analytic in $\mathbb{D}$ (by the
usual comparison test) and hence they are members of $\mathcal{A}$.
However, they need not be univalent. For example, the function
\[f(z)=z-2z^2-3z^3-4z^4-\cdots= 2z-\frac{z}{(1-z)^2}\] satisfies the
inequality $|a_n|\leq n$ but its derivative vanishes inside
$\mathbb{D}$ and therefore the function $f$ is not univalent in
$\mathbb{D}$. In 1970, Gavrilov \cite{gavri}  showed that  the
radius of univalence of functions satisfying the inequality
$|a_n|\leq n$ is the real root of the equation $2(1-r)^3-(1+r)=0$
while, for the functions whose coefficients satisfy $|a_n|\leq M$,
the radius of univalence is $1-\sqrt{M/(1+M)}$. Later, in 1982,
Yamashita showed that the radius of univalence obtained by Gavrilov
is also the same as the radius of starlikeness of the corresponding
functions. He also found lower bounds for the radii of convexity for
these functions. Recently, in 2006, Graham \emph{et al.}\
\cite[Theorem 4.2 and Lemma 5.6]{graham} considered the
corresponding radius problems for holomorphic mappings on the unit
ball in~$\mathbb{C}^n$.  Kalaj,   Ponnusamy, and  Vuorinen
\cite{ponnu} have investigated related problems for harmonic
functions.  In this paper, several related radius problems for the
following classes of functions will be investigated.

For $ 0\leq \alpha<1$, let $\mathcal{S}^* (\alpha)$ and
$\mathcal{C}(\alpha)$ be subclasses of $\mathcal{S}$ consisting of
starlike functions of order $\alpha$ and convex functions of order
$\alpha$, respectively defined analytically by the following
equalities:
\[\mathcal{S}^*(\alpha):=\left\{ f\in\mathcal{S}: \RE \left( \frac{zf'(z)}{f(z)} \right)
> \alpha\right\}, \
\text{ and } \ \mathcal{C}(\alpha):=\left\{f\in\mathcal{S}: \RE
\left(1+ \frac{zf''(z)}{f'(z)} \right) >\alpha \right\}.
\] The classes $\mathcal{S}^*:=\mathcal{S}^*(0)$ and
$\mathcal{C}:=\mathcal{C}(0)$ are the familiar classes of starlike
and convex functions respectively. Closely related are the following
classes of functions:
\[\mathcal{S}^*_\alpha  :=\left\{ f\in\mathcal{S}:   \left| \frac{zf'(z)}{f(z)} -1
\right| < 1- \alpha\right\}, \quad \text{ and } \quad
 \mathcal{C}_\alpha :=\left\{ f\in\mathcal{S}:  \left|\frac{zf''(z)}{f'(z)}
\right| < 1-\alpha\right\} .
\] Note that $ \mathcal{S}^*_\alpha \subseteq\mathcal{S}^*(\alpha)$
and $\mathcal{C}_\alpha\subseteq \mathcal{C}(\alpha)$.

A function $f\in \mathcal{S}$ is \emph{uniformly convex} if $f$ maps
every circular arc $\gamma$ contained in $\mathbb{D}$ with center
$\zeta\in\mathbb{D}$  onto a convex arc. The class of all uniformly
convex functions, introduced by Goodman \cite{goodman}, is denoted
by $\mathcal{UCV}$. R\o nning \cite[Theorem 1, p.\ 190]{MR1128729},
and  Ma and Minda \cite[Theorem 2, p.\ 162]{MR1182182},
independently  showed that  $ f \in  \mathcal{S} $ is uniformly
convex if and only if
\begin{equation*}
\RE \left( 1 + \frac{zf''(z)}{f'(z)} \right) >
\left|\frac{zf''(z)}{f'(z)} \right| \quad ( z \in  \mathbb{D}) .
\end{equation*}
R\o nning  \cite{MR1128729} also considered the class
$\mathcal{S}_P$ of parabolic starlike functions consisting of
functions $f \in \mathcal{A} $ satisfying
\[   \label{eqsp} \RE \left(    \frac{zf'(z)}{f(z)}     \right)     > \left|
\frac{zf'(z)}{f(z)} - 1\right|  \quad(  z\in \mathbb{D}) . \] In
other words, the class $\mathcal{S}_P$ consists of function $f=zF'$
where $ F \in \mathcal{UCV}$. For a recent survey on uniformly
convex functions, see \cite{ravi}.


For a fixed $b$ with $0\leq b\leq 1$, let $\mathcal{A}_b$ denote the
class of all analytic functions $f$   of the form
\[f(z)=z+a_2z^2+a_3z^3+a_4z^4+\cdots\quad (|a_2|=2b,\ z\in \mathbb{D}).\]
The second coefficient of univalent functions determines important
properties such as growth and distortion estimates.  For recent
investigation of functions with fixed second coefficients, see
\cite{vravi2,vravi}. For $0\leq \alpha <1$, the sharp radii of
starlikeness and convexity of order $\alpha$  are obtained for
functions $f\in\mathcal{A}_b$ satisfying the condition  $|a_n|\leq n
$ or $|a_n|\leq M$ ($M>0$) for $n\geq 3$. Special case ($\alpha=0$)
of the results  shows that the lower bounds for the radii of
convexity obtained by Yamashita \cite{yama} are indeed sharp. The
coefficient inequalities are natural in the sense that the
inequality $|a_n|\leq n $ is satisfied by univalent functions and
while the inequality $|a_n|\leq M$ is satisfied by functions which
are bounded by $M$. For a function $p(z)=1+c_1z+c_2z^2+\cdots$ with
positive real part, it is well-known that $|c_n|\leq 2$ and so if
$f\in\mathcal{A}$ and $\RE f'(z)>0$, then $|a_n|\leq 2/n$. In view
of this, the determination of the radius of starlikeness and the
radius of convexity of functions whose coefficients satisfy  the
inequality $|a_n|\leq M/n$ is also investigated.  A corresponding
radius problem for certain function $p(z)=1+c_1z+c_2z^2+\cdots$ with
coefficients satisfying the conditions $|c_1|=2b$, $0\leq b\leq 1$
and $|c_n|\leq 2M$ ($M>0$) is also investigated.

\section{Radii of starlikeness of order \boldmath $\alpha$\unboldmath\ and parabolic starlikeness}
In this section,  the sharp $\mathcal{S}^*(\alpha)$-radius and the
sharp $\mathcal{S}^*_\alpha$-radius for  $0\leq \alpha <1$ as well
as the sharp $\mathcal{S}_P$-radius  are obtained for functions
$f\in\mathcal{A}_b$ satisfying one of the conditions  $|a_n|\leq n
$, $|a_n|\leq M$ or $|a_n|\leq M/n$ ($M>0$) for $n\geq 3$.

\begin{theorem}\label{th21}
Let $f\in \mathcal{A}_b$ and $|a_n|\leq n$ for $n\geq 3$. Then $f$
satisfies the inequality
\begin{equation}\label{th1eq1}
\left|\frac{zf'(z)}{f(z)}-1\right|\leq 1-\alpha \quad (|z|\leq r_0)
\end{equation}
where $r_0=r_0(\alpha)$ is the real  root in $(0,1)$ of the equation
\begin{equation}\label{th1e2}
1-\alpha+(1+\alpha)r=2\big(  1-\alpha +(2-\alpha)(1-b)r\big)(1-r)^3.
\end{equation} The number $r_0(\alpha)$ is also the
radius of starlikeness of order $\alpha$. The number $r_0(1/2)$ is
the radius of parabolic starlikeness of the given functions.  The
results are all sharp.
\end{theorem}


\begin{proof}
If
\begin{equation}\label{th1eq2}
\sum_{n=2}^\infty (n-\alpha)|a_n|r_0^{n-1} \leq 1-\alpha,
\end{equation}
 then
 the function $f(z)=z+\sum_{n=2}^\infty
a_nz^n$ satisfies, on $|z|=r_0$,
\begin{align*}
|zf'(z)-f(z)|-(1-\alpha)|f(z)| &\leq \sum_{n=2}^\infty
(n-1)|a_n||z|^n -(1-\alpha)(|z|-\sum_{n=2}^\infty |a_n|\,|z|^n )\\
& = -(1-\alpha)|z|+\sum_{n=2}^\infty (n-\alpha)|a_n||z|^n \\
& \leq r_0 \left( -(1-\alpha) +\sum_{n=2}^\infty
(n-\alpha)|a_n|r_0^{n-1}\right)\\
&\leq 0.
\end{align*}
This shows that the condition \eqref{th1eq2} is a sufficient
condition for the inequality \eqref{th1eq1} to hold.  Using
$|a_2|=2b$ for the function $f\in \mathcal{A}_b$, and the inequality
$|a_n|\leq n$ for $n\geq 3$, it follows that, for $|z|\leq r_0$,
\begin{align*}
\sum_{n=2}^\infty (n-\alpha)|a_n|\,|z|^{n-1} & \leq \sum_{n=2}^\infty (n-\alpha)|a_n|r_0^{n-1} \\
& \leq 2(2-\alpha)br_0+\sum_{n=3}^\infty
n^2r_0^{n-1}-\alpha\sum_{n=3}^\infty
nr_0^{n-1} \\
&=2(2-\alpha)br_0+\frac{1+r_0}{(1-r_0)^3}-1-4r_0-\alpha\left(\frac{1}{(1-r_0)^2}-1-2r_0\right)\\
&=\alpha-1-2(2-\alpha)(1-b)r_0+\frac{(1+r_0)-\alpha(1-r_0)}{(1-r_0)^3}\\
&= 1-\alpha
\end{align*}
provided $r_0$ is the root of the  Equation~\eqref{th1e2} in the
hypothesis of the theorem. The  Equation~\eqref{th1e2} has clearly a
root in (0,1).

The function $f_0$ given by
\begin{equation}\label{sharpf0}
f_0(z)=2z+2(1-b)z^2-\frac{z}{(1-z)^2}=z-2bz^2-3z^3-4z^4-\cdots
\end{equation}
satisfies the hypothesis of the theorem and, for this function, we
have
\[\frac{zf_0'(z)}{f_0(z)}-1
=\frac{2(1-b)z(1-z)^3-2z}{(2+2(1-b)z)(1-z)^3-(1-z)}. \] For $z=r_0$,
we have
\begin{equation}\label{th1e4}
 \left|\frac{zf_0'(z)}{f_0(z)}-1\right|
=\frac{2r_0-2(1-b)r_0(1-r_0)^3}{(2+2(1-b)r_0)(1-r_0)^3-(1-r_0)}=1-\alpha.
\end{equation}
This shows that the radius $r_0$ of functions to satisfy
\eqref{th1eq1} is sharp. The numerator of the rational function in
the middle of \eqref{th1e4} is positive as $0\leq 1-b \leq 1$ and
$0\leq (1-r_0)<1$ which shows that $(1-b)(1-r_0)^3<1$. The
denominator expression is also positive as $(2+2(1-b)r_0)(1-r_0)^2
\geq 2(1-r_0)^2 >1$. The inequality $2(1-r_0)^2>1$ is in fact
equivalent to $r_0 < 1-1/\sqrt{2}=0.292893$. This inequality holds
as $r_0=r_0(\alpha) \leq r_0(0)=0.1648776 $.

Since the functions satisfying \eqref{th1eq1} are starlike of order
$\alpha$, the  radius of starlikeness is at least $r_0(\alpha)$.
However, this radius is also sharp for the same function $f_0$ as
\begin{equation}\label{th1e5}
\RE\left(\frac{zf_0'(z)}{f_0(z)}\right) = \alpha \quad (z=r_0).
\end{equation}

The inequality
\[ \left|\frac{zf'(z)}{f(z)}-1\right|\leq \frac{1}{2}\]
is sufficient (see \cite{ravi}) for the function  to be parabolic
starlike and hence the radius of parabolic starlikeness is at least
$r_0(1/2)$.
 The  Equations~\eqref{th1e4} and
\eqref{th1e5} with $\alpha=1/2$ shows that
\begin{equation}\label{th1e6}
 \left|\frac{zf_0'(z)}{f_0(z)}-1\right| =\frac{1}{2}=
 \RE\left(\frac{zf_0'(z)}{f_0(z)}\right)   \quad (z=r_0),
\end{equation}
and hence  the radius of parabolic starlikeness is sharp.
\end{proof}

\begin{corollary}\label{cor22}
The radius of starlikeness of order $\alpha$ of functions whose
coefficients satisfy $|a_n|\leq n$ for all $n\geq 2$ is the real
root in $(0,1)$ of the equation
\[2(1-\alpha)(1-r)^3=1-\alpha+(1+\alpha)r.\] In particular, the
radius of starlikeness is given by \[ r_0(0)= 1+\frac{1}{6^{2/3}}
\left( (\sqrt{330}-18)^{1/3}- (\sqrt{330}+18 )^{1/3} \right) \approx
0.164878.
\]  The radius of starlikeness of order 1/2 is the same as  the radius of
parabolic starlikeness and it is given by
\[ r_0(1/2)=
1+\frac{1}{\sqrt{2}}\left( \left(3-2\sqrt{2}
\right)^{1/3}-\left(3+2\sqrt{2}\right)^{1/3} \right) \approx
0.120385.\]  The results are sharp.
\end{corollary}

\begin{corollary}\label{cor23}
The radius of starlikeness of order $\alpha$ of functions whose
coefficients satisfy $a_2=0$ and $|a_n|\leq n$ for all $n\geq 3$ is
the real root in $(0,1)$  of the equation
\[2(1-\alpha+(2-\alpha)r)(1-r)^3=1-\alpha+(1+\alpha)r.\] In particular,
the radius of starlikeness is the root $r_0\approx 0.253571$    of
the equation \[ 2 (1 + 2 r) (1 - r)^3 = 1 + r .\] The radius of
starlikeness of order 1/2 which is  the same as  the radius of
parabolic starlikeness is $r_0=1-\sqrt[3]{1/2}\approx 0.206299$. The
results are sharp.
\end{corollary}

\begin{remark}
It is clear from Corollaries  \ref{cor22} and \ref{cor23} that the
various radii are improved if the second coefficient of the function
vanishes.
\end{remark}


%

\begin{theorem} \label{th22}
Let $f\in \mathcal{A}_b$ and $|a_n|\leq M$ for $n\geq 3$. Then $f$
satisfies the condition \eqref{th1eq1}  where $r_0=r_0(\alpha)$ is
the real  root in $(0,1)$ of the equation
\[ M(1-\alpha+\alpha r)=\big((1+M)(1-\alpha)-(2-\alpha)(2b-M)r\big)(1-r)^2\]
The number $r_0(\alpha)$ is also the radius of starlikeness of
order $\alpha$. The number $r_0(1/2)$ is the radius of parabolic
starlikeness of the given functions.   The results are all sharp.
\end{theorem}

\begin{proof}
Using   $|a_2|=2b$ for the function $f\in \mathcal{A}_b$, and the
inequality $|a_n|\leq M$ for $n\geq 3$, a calculation shows that,
for $|z|\leq r_0$,
\begin{align*}
\sum_{n=2}^\infty (n-\alpha)|a_n|\,|z|^{n-1} & \leq
\sum_{n=2}^\infty (n-\alpha)|a_n|r_0^{n-1}\\
& \leq 2(2-\alpha)br_0+M\left(\sum_{n=3}^\infty
nr_0^{n-1}-\alpha\sum_{n=3}^\infty
r_0^{n-1} \right)\\
&= 2(2-\alpha)br_0+M\left( \frac{1}{(1-r_0)^2}-1-2r_0-\alpha\left(\frac{
1}{1-r_0} -1-r_0 \right)\right)\\
&= (2-\alpha)(2b-M)r_0-M(1-\alpha)+M\frac{1-\alpha+\alpha r_0}{(1-r_0)^2}\\
&=1-\alpha
\end{align*}
where $r_0$ is as stated in the hypothesis of the theorem. Thus, the
function $f$ satisfies the condition \eqref{th1eq1}. The other two
results follow easily.

The results are sharp for the function $f_0$ given by
\begin{equation}\label{th2e3}
f_0(z)=z-2bz^2-M(z^3+z^4+\cdots)= z-2bz^2-\frac{Mz^3}{1-z}.
\end{equation} A calculation shows that
\[ \frac{zf_0'(z)}{f_0(z)}-1 = -
\frac{ 2bz+\frac{2Mz^2}{1-z}+\frac{Mz^3}{(1-z)^2}}{
1-2bz-\frac{Mz^2}{1-z}}
\]
At the point $z=r_0$, the function $f_0$ satisfies
\[ \RE\left(\frac{zf_0'(z)}{f_0(z)}\right)= 1-
\frac{ 2br_0+\frac{2Mr_0^2}{1-r_0}+\frac{Mr_0^3}{(1-r_0)^2}}{
1-2br_0-\frac{Mr_0^2}{1-r_0}} =\alpha.
\]
Since $\alpha<1$, the last equation shows that  the denominator of
the rational expression in the middle is positive. This leads to the
following equality:
\[ \left|\frac{zf_0'(z)}{f_0(z)}-1\right| =
\frac{ 2br_0+\frac{2Mr_0^2}{1-r_0}+\frac{Mr_0^3}{(1-r_0)^2}}{
1-2br_0-\frac{Mr_0^2}{1-r_0}} =1-\alpha. \] Also the Equation
  \eqref{th1e6} holds. This proves the sharpness of the
results.
\end{proof}

\begin{corollary}Let $f\in \mathcal{A}$ and $|a_n|\leq M$ for $n\geq
2$. Then $f$ satisfies the condition \eqref{th1eq1}  where
$r_0(\alpha)$ is the real  root in $(0,1)$ of the equation
\[ M(1-\alpha+\alpha r)=(1+M)(1-\alpha)(1-r)^2. \]
The number $r_0(\alpha)$ is also the radius of starlikeness of order
$\alpha$. The number $r_0(1/2)$ is the radius of parabolic
starlikeness of the given functions.   The results are all sharp.
\end{corollary}

\begin{remark}
The radius of starlikeness of the functions $f$ with $|a_n|\leq M$
given by $r_0=1-\sqrt{M/(1+M)}$  is the  root in (0,1) of the
equation
\[ M=(1+M)(1-r)^2.\] When the second coefficient $a_2=0$, the radius
of starlikeness $r_1$ is the root in (0,1) of the equation \[
M=(1+M+2Mr)(1-r)^2.\] Clearly, $r_1>r_0$.
\end{remark}

\begin{theorem}\label{th23}
Let $f\in \mathcal{A}_b$ and $|a_n|\leq M/n$ for $n\geq 3$. Then $f$
satisfies the condition \eqref{th1eq1}  where $r_0=r_0(\alpha)$ is
the real  root in $(0,1)$ of the equation
\[ 2M(1+\alpha(1-r)(\log(1-r))/r)=(2(1+M)(1-\alpha)+(2-\alpha)(M-4b)r)(1-r)\]
The number $r_0(\alpha)$ is also the radius of starlikeness of order
$\alpha$. The number $r_0(1/2)$ is the radius of parabolic
starlikeness of the given functions.   The results are all sharp for
the function $f_0$ given by \[ f_0(z):= (1+M)z+(M/2-2b)z^2+M
\log(1-z).\]
\end{theorem}

The logarithm in the above equation is the branch that takes the
value 1 at $z=0$.  Proof of this theorem is   omitted as it is
similar to those of Theorems~\ref{th21} and \ref{th22}.

\section{Radii of convexity and  uniform convexity}
In this section,  the sharp $\mathcal{C}(\alpha)$-radius and the
sharp $\mathcal{C}_\alpha$-radius for  $0\leq \alpha <1$ as well as
the sharp $\mathcal{UCV}$-radius for functions $f\in\mathcal{A}_b$
satisfying the condition $|a_n|\leq n $ or $|a_n|\leq M$ ($M>0$) for
$n\geq 3$ are obtained.

\begin{theorem}Let $f\in \mathcal{A}_b$ and $|a_n|\leq n$ for $n\geq
3$. Then $f$ satisfies the condition
\begin{equation}\label{th3eq1}
\left|\frac{zf''(z)}{f'(z)} \right|\leq 1-\alpha \quad (|z|\leq r_0)
\end{equation}
where $r_0=r_0(\alpha)$ is the real  root in $(0,1)$ of the equation
\begin{equation}\label{th3e2}
2\big(1-\alpha+2(2-\alpha)(1-b)r\big)(1-r)^4=1-\alpha+4r+(1+\alpha)r^2
\end{equation} The number $r_0(\alpha)$ is
also the radius of convexity of order $\alpha$. The number
$r_0(1/2)$ is the radius of uniform convexity of the given
functions. The results are all sharp.
\end{theorem}



\begin{proof} A function $f$ satisfies \eqref{th3eq1} if and only if
$zf'$ satisfies \eqref{th1eq1}. In view of this and the inequality
\eqref{th1eq2}, the inequality
\begin{equation}\label{th3eq2}
\sum_{n=2}^\infty n (n-\alpha)|a_n||z|^{n-1} \leq 1-\alpha, \quad
(|z|\leq r_0)
\end{equation}
is sufficient for function $f$ to satisfy \eqref{th3eq1}. Let $r_0$
be the root in $(0,1)$ of the Equation \eqref{th3e2}. Now, for
$|z|\leq r_0$,
\begin{align*}
\sum_{n=2}^\infty n (n-\alpha)|a_n||z|^{n-1} & \leq
\sum_{n=2}^\infty n (n-\alpha)|a_n|r_0^{n-1}\\
 & \leq  4 (2 - \alpha) b r_0+\sum_{n=3}^\infty (n - \alpha) n^2 r_0^{n - 1} \\
&  = 4 (2 - \alpha) b r_0  +\left(\frac{1+4r_0+r_0^2}{(1-r_0)^4}-1-8r_0\right)
-\alpha\left(\frac{1+r_0}{(1-r_0)^3}-1-4r_0\right)\\
& = -\big(1-\alpha+4(2-\alpha)(1-b)r_0\big)+\frac{1-\alpha+4r_0+(1+\alpha)r_0^2}{(1-r_0)^4}\\
& = 1-\alpha.
 \end{align*}

To prove the sharpness, consider the function $f_0$ defined by \eqref{sharpf0}. For
this function, a calculation shows that
\[ \frac{zf_0''(z)}{f_0'(z)}=
 \frac{ 4 (1-b)z-\frac{4z}{(1-z)^3}-\frac{6 z^2}{(1-z)^4} }
 {2+4 (1-b) z-\frac{1}{(1-z)^2}-\frac{2 z}{(1-z)^3}}
.\] If $r_0$ is the root of the equation  \eqref{th3e2}, then, at the point $z=r_0$,
\[ \RE\left(\frac{zf_0''(z)}{f_0'(z)}\right)=
 \frac{ 4 (1-b)r_0-\frac{4r_0}{(1-r_0)^3}-\frac{6 r_0^2}{(1-r_0)^4} }
 {2+4 (1-b) r_0-\frac{1}{(1-r_0)^2}-\frac{2 r_0}{(1-r_0)^3}}=\alpha-1
.\]The denominator of the rational function in the middle of the
equation above is positive while the numerator is negative. Noting
this, it also follows that, at the point $z=r_0$,
\[\left|\frac{zf_0''(z)}{f_0'(z)}\right| =
 \frac{ -4 (1-b)r_0+\frac{4r_0}{(1-r_0)^3}+\frac{6 r_0^2}{(1-r_0)^4} }
 {2+4 (1-b) r_0-\frac{1}{(1-r_0)^2}-\frac{2 r_0}{(1-r_0)^3}}=1-\alpha.\]
 In the case of $\alpha=1/2$, the equation \eqref{th1e6} also holds.
\end{proof}

The special case where $b=1$ is important and it is stated as a corollary below.

\begin{corollary}\label{cor33}
Let $f\in \mathcal{A}$ and $|a_n|\leq n$ for $n\geq 2$. Then $f$
satisfies the condition \eqref{th3eq1} where $r_0=r_0(\alpha)$ is
the real  root in $(0,1)$ of the equation
\begin{equation}\label{th3e2c2}
2(1-\alpha)(1-r)^4=1-\alpha+4r+(1+\alpha)r^2
\end{equation} The number $r_0(\alpha)$ is
also the radius of convexity of order $\alpha$. The number
$r_0(1/2) \approx 0.064723$ is the radius of uniform convexity of the given
functions. The results are all sharp.
\end{corollary}

\begin{remark} For $\alpha=0$, the Equation \eqref{th3e2c2}
reduces to \[ 2 (1 - r)^4 = (1 + 4 r + r^2). \] The root of this
equation in (0,1) is approximately  0.09033. Our result shows that
radius of convexity obtained by Yamashita \cite[Theorem 2]{yama} is
sharp.
\end{remark}

\begin{corollary}\label{cor34}
Let $f\in \mathcal{A}$, $a_2=0$ and $|a_n|\leq n$ for $n\geq 3$.
Then $f$ satisfies the condition \eqref{th3eq1} holds where
$r_0=r_0(\alpha)$ is the real  root in $(0,1)$ of the equation
\begin{equation}\label{th3e2c}
2(1-\alpha+2(2-\alpha)r)(1-r)^4=1-\alpha+4r+(1+\alpha)r^2
\end{equation} The number $r_0(\alpha)$ is
also the radius of convexity of order $\alpha$. The number $r_0(1/2)
\approx 0.125429$ is the radius of uniform convexity of the given
functions. The results are all sharp.
\end{corollary}

\begin{remark}It is easy to see from Corollaries~\ref{cor33} and
\ref{cor34} that the radius of convexity of order $\alpha$ improves
when $a_2=0$. In the particular case $\alpha=0$, the root of the
Equation \eqref{th3e2c2} is $r_0(0) \approx 0.0903331$ while the
Equation \eqref{th3e2c} has the root  $r_0(0) \approx 0.155972$.
\end{remark}

\begin{theorem}Let $f\in \mathcal{A}_b$ and $|a_n|\leq M$ for $n\geq
3$. Then $f$ satisfies the condition \eqref{th3eq1}  where $r_0=r_0(\alpha)$ is
the real  root in $(0,1)$ of the equation
\[ \big((1- \alpha)(1+M)-2(2- \alpha)(2b-M)r\big)(1-r)^3=M\big(1- \alpha+(1+ \alpha)r\big).
\] The number $r_0(\alpha)$ is also the radius of convexity of
order $\alpha$. The number $r_0(1/2)$ is the radius of uniform convexity
of the given functions.   The results are all sharp.
\end{theorem}

\begin{proof}
Using   $|a_2|=2b$ for the function $f\in \mathcal{A}_b$, and the
inequality $|a_n|\leq M$ for $n\geq 3$, a calculation shows that,
for $|z|\leq r_0$,
\begin{align*}
\sum_{n=2}^\infty n(n-\alpha)|a_n|\,|z|^{n-1} & \leq
\sum_{n=2}^\infty n (n-\alpha)|a_n|r_0^{n-1}\\
& \leq 4(2-\alpha)br_0+M\left(\sum_{n=3}^\infty
n^2r_0^{n-1}-\alpha\sum_{n=3}^\infty
n r_0^{n-1} \right)\\
&= 4(2-\alpha)br_0+M\left( \frac{1+r_0}{(1-r_0)^3}-1-4r_0-\alpha\left( \frac{
1}{(1-r_0)^2} -1 -2r_0\right) \right)\\
& = -M(1- \alpha)+2(2- \alpha)(2b-M)r_0+M\left(\frac{1-\alpha+(1+\alpha)r_0}{(1-r_0)^3}\right)\\
&=1-\alpha
\end{align*}
where $r_0$ is as stated in the hypothesis of the theorem. Thus, the
function $f$ satisfies the condition \eqref{th3eq1}. The other two
results follow easily. The results are sharp for the function $f_0$
given by \eqref{th2e3}.
\end{proof}

\begin{corollary}Let $f\in \mathcal{A}$ and $|a_n|\leq M$ for $n\geq
2$. Then $f$ satisfies the condition \eqref{th3eq1}  where $r_0=r_0(\alpha)$ is
the real  root in $(0,1)$ of the equation
\[ (1- \alpha)(1+M) (1-r)^3=M\big(1- \alpha+(1+ \alpha)r\big).
\] The number $r_0(\alpha)$ is also the radius of convexity of
order $\alpha$. The number $r_0(1/2)$ is the radius of uniform convexity
of the given functions.   The results are all sharp.
\end{corollary}

\begin{remark} For $\alpha=0$, the Equation \eqref{th3e2}
reduces to \[  (1+M^{-1}) (1-r)^3= 1+r .
\]   Our result again shows that radius of convexity obtained
by Yamashita \cite[Theorem 2]{yama} is sharp.
\end{remark}

\begin{remark} The problem of determining the radius of convexity of
functions satisfying $|a_n|\leq M/n$ is the same as the
determination of radius of starlikeness of functions satisfying the
inequality $|a_n|\leq M$. The latter problem is investigated in
Theorem~\ref{th23}.\end{remark}

\section{Carath\'eodory functions }
An analytic function $p$ of the form $p(z)=1+c_1z+c_2z^2+\cdots$ is
called a Carath\'eodory function if $\RE p(z)>0$ for all $z\in
\mathbb{D}$. The class of all such functions is denoted by
$\mathcal{P}$.  For such functions $p\in \mathcal{P}$, it is
well-known that $|c_n|\leq 2$. Denote the class of all
Carath\'eodory functions satisfying  the inequality $\RE
p(z)>\alpha$ for some $0\leq \alpha <1$ by $\mathcal{P}(\alpha)$. It
is  easy to see that $|c_n|\leq 2(1-\alpha)$ for
$p\in\mathcal{P}(\alpha)$. In this section, we determine
$\mathcal{P}(\alpha)$-radius of functions satisfying the inequality
$|c_n|\leq 2M $ for $n\geq 3$ with $|c_2|=2b$  fixed. The proof of
the  following result  is straightforward and the details are
omitted.

\begin{theorem}Let    $p$ be  an analytic function of the form
$p(z)=1+c_1z+c_2z^2+\cdots$ with   $|c_2|=2b$ and $|c_n|\leq 2M $
for $n\geq 3$. Then \[ |p(z)-1|\leq 1-\alpha \quad (|z|\leq r_0)\]
where \[ r_0=r_0(\alpha) =
\frac{2(1-\alpha)}{1-\alpha+2b+\sqrt{(1-\alpha+2b)^2+8(1-\alpha)(M-b)}}.
\] Also $\RE p(z)>\alpha$ for $|z|\leq r_0(\alpha)$. These results
are sharp for the function $p_0$ given by
\[ p_0(z)=1-2bz-2M\frac{z^2}{1-z}. \]
\end{theorem}

\end{document}